\documentclass[12pt]{article}
\usepackage{amsmath}
\usepackage{latexsym,amssymb,eucal}
\newtheorem{thm}{Theorem}[section]
\newtheorem{lem}[thm]{Lemma}

\newtheorem{exam}[thm]{Example}
\newtheorem{df}[thm]{Definition}
\newtheorem{cor}[thm]{Corollary}
\newcommand{\id}{\mathrm{id}}
\newcommand{\Ad}{\mathrm{Ad}\,}

\newcommand{\tM}{\tilde{\mathcal{M}}}
\newcommand{\tN}{\bar{\mathcal{N}}}

\newcommand{\tQ}{\bar{\mathcal{Q}}}
\newcommand{\cA}{\mathcal{A}}
\newcommand{\cB}{\mathcal{B}}
\newcommand{\cM}{\mathcal{M}}
\newcommand{\cP}{\mathcal{P}}
\newcommand{\cQ}{\mathcal{Q}}
\newcommand{\cR}{\mathcal{R}}
\newcommand{\cS}{\mathcal{S}}
\newcommand{\cT}{\mathcal{T}}
\newcommand{\cN}{\mathcal{N}}

\newcommand{\Aut}{\mathrm{Aut}}

\newcommand{\End}{\mathrm{End}}

\newcommand{\Mor}{\mathrm{Mor}}

\begin{document}

\title{Classification of Roberts actions of strongly amenable
C$^*$-tensor categories on the injective factor of type III$_1$} 
\author{Toshihiko MASUDA\footnote{Supported by 
JSPS KAKENHI Grant Number 16K05180.} \\
 Graduate School of Mathematics, Kyushu University \\
 744, Motooka, Nishi-ku, 
Fukuoka, 819-0395, JAPAN \\
e-mail address: masuda@math.kyushu-u.ac.jp}
\date{}

\maketitle

\begin{abstract}
In this paper, we generalize Izumi's result on uniqueness of realization
 of finite C$^*$-tensor categories in the endomorphism category of the
 injective factor of type III$_1$
for finitely generated strongly amenable C$^*$-tensor categories by
 applying Popa's classification theorem of strongly amenable subfactors
 of type III$_1$.
\end{abstract}

\section{Introduction}\label{sec:intro}

Since V. F. R. Jones initiated the theory of index for subfactors
\cite{J-ind}, the theory of subfactors has been developed surprisingly by
involving many area of mathematics, e.g., knot theory, low dimensional
topology, mathematical physics.
One of important points is the 
connection between subfactor theory and tensor category theory.
If a subfactor is given, then a C$^*$-tensor
category naturally arises as a bimodule category in type II$_1$ setting,
and endomorphism category in type III setting. These C$^*$-tensor
categories contain rich informations about the structure of subfactors.
See \cite{EK-book} for details of relation between subfactor theory and tensor categories.

Conversely 
any amenable
C$^*$-tensor category can be realized as a full subcategory of a
bimodule category of an AFD type
II$_1$ factor by the result of Hayashi-Yamagami \cite{HaYa-amenable}. 
(Even if one drops the assumption on amenability, one can
still realize it as a bimodule category of some non-AFD type II$_1$ factor \cite{Yama-real}.)
By tensor product trick, one can realize an amenable
C$^*$-tensor category as a full subcategory of endomorphisms of an
injective type III$_1$ factor. So natural question is the uniqueness of
realization, or more generally
how C$^*$-tensor functors between two C$^*$-tensor categories are realized in the setting of
operator algebras.

In \cite[Theorem 2.2]{Izumi-neargrp}, Izumi showed that the equivalence
of finite C$^*$-tensor categories is given by an isomorphism of
injective factors of type III$_1$, which implies the essential
uniqueness of embedding of a given finite C$^*$-tensor category into a
C$^*$-tensor category of endomorphisms of the injective factor of type III$_1$. 
This result can be
regarded  as the  classification of ``actions'' of C$^*$-tensor categories. For
example, if a C$^*$-tensor category comes from  a finite group, then his result means
the  uniqueness of free finite group actions on the injective factor of
type III$_1$.  

Izumi's proof is based on Popa's deep classification result of amenable
subfactors \cite{Po-amen}, \cite{Po-tani}, \cite{Po-CBMS}. It is
well-known that classification of subfactors yields that of discrete
amenable group actions on injective factors by considering locally
trivial subfactors \cite[\S 5.1.5]{Po-amen}.
However, one must look
carefully on behavior on standard invariants to obtain sufficient
classification result. In fact, one get only outer conjugacy of actions
from an isomorphism of two locally trivial subfactors. 
In \cite[Theorem 2.2]{Izumi-neargrp}, Izumi uses variants of
locally trivial subfactors for finite C$^*$-tensor categories.
To obtain sufficient results, he chose an isomorphism of subfactors which
preserves a given tunnel with finite length, and investigated standard
invariants carefully. 
However his method is valid
only for finite C$^*$-tensor categories. 

In this paper, we generalize Izumi's result, and show a similar result
for finitely generated strongly amenable C$^*$-tensor categories. Our
main tools are the notion of Loi invariant \cite{Loi-auto} and
technique used in \cite{M-ExtIJM}. By applying Popa's classification of
subfactor of type III$_1$ and investigating an isomorphism of standard
invariants carefully, we deduce the main result.

\textbf{Aknowledgement.} The author is grateful to Professor Reiji
Tomatsu for fruitful discussion and suggesting him for simplification of proof.

\section{Loi invariant}\label{sec:half}

Our standard references are \cite{EK-book} for subfactor theory,
\cite{Iz-fusion} for sector theory and \cite{Nes-Tu} for C$^*$-tensor categories.

In this section, we recall basic facts on the Loi invariant \cite{Loi-auto}.

Let $\cN\subset \cM$ and $\cQ\subset \cP$ be isomorphic subfactors with
finite index, and $\alpha\colon\cN\subset \cM\rightarrow \cQ\subset \cP$ be
an isomorphism. 
Fix tunnels 
\[
\cM\supset \cN_1\supset \cN_2\supset \cN_3\supset \cdots, \,\,\,
\cP\supset \cQ_1\supset \cQ_2\supset \cQ_3\supset \cdots. 
    \]
Let $e_k^{\cN}\in \cN_k$ (resp. $e_k^{\cQ}\in \cQ_k$)
be a Jones projection for $\cN_{k+2}\subset \cN_{k+1}$ 
(resp. $\cQ_{k+2}\subset \cQ_{k+1}$).
Here $\cN_0=\cM$, $\cQ_0=\cP$.
Then for any $k\geq 2$, we can choose a unitary $u\in \cQ$ such that 
$\Ad u\circ \alpha(e_l^{\cN})=e_l^{\cQ}$, $0\leq l\leq k-2$. Hence $\Ad
u\circ \alpha$ preserves a tunnel 
\[
\cM\supset \cN_1\supset \cN_2\supset \cN_3\cdots\supset \cN_k,
  \]
i.e., $\Ad u\circ \alpha(\cN_l)=\cP_l$, $0\leq l \leq k$. Thus $\Ad
u\circ \alpha $ induces an isomorphism $\cN_l'\cap \cM\rightarrow
\cQ_l'\cap \cP$, $0\leq l\leq k$. It is shown that $\Ad u\circ \alpha
|_{\cN'_l\cap \cM}$ is independent from the choice of $u$. Thus
$\alpha_{l}:=\Ad u\circ\alpha |_{\cN'_l\cap \cM}$ is well-defined, i.e.,
$\alpha_l$ does not depend on $u$ nor $k$,
and
the following definition is justified.
\begin{df}[{\cite[p.286]{Loi-auto}}]
 The Loi invariant $\Phi(\alpha)$ of $\alpha$ is defined by 
$\Phi(\alpha)=\{\alpha_k\}_{k=0}^\infty$.
\end{df}
We clarify the behavior of $\Phi(\alpha)$ by choice of tunnels.
Let
\[
\cM\supset \cN_1\supset \tN_2\supset \tN_3\cdots,
\cP\supset \cQ_1\supset \tQ_2\supset \tQ_3\cdots. 
    \]
be different choice of tunnels. Denote Jones projections for these
tunnels by $e_k^{\tN}\in \tN_k$, $e_k^{\tQ}\in \tQ_k$.
For any $k\geq 2$, there exist
unitaries $v\in \cN$,  $w\in \cQ$ such that  $ve_{l}^{\cN}v^*=e_l^{\tN}$,
$we_{l}^{\tQ}w^*=e_l^{\tQ}$, $0\leq l\leq k-2$. This implies $v\cN_l
v^*=\tN_l$, $w\cQ_lw^*=\tQ_l$, $0\leq l\leq k$.
By this unitary
perturbation, we get a Loi invariant $\{\Ad w \circ \alpha_k\circ \Ad v^*\}_{k=1}^\infty$
for this new tunnel, and this does not depend on the choice of unitaries
$v,w$.
In this sense, the Loi invariant $\Phi(\alpha)$ is defined canonically
for an isomorphism $\alpha\colon\cN\subset 
\cM\rightarrow \cQ\subset \cP$,  and it defines 
an isomorphism of standard invariants. (Later we
will explain details of an isomorphism of standard invariants.)

From now on, we assume that involved factors are of type III. 
Let $\End_0(\cM)$ be a set of unital endomorphisms of $\cM$ with finite
statistical dimension.
For $\pi_1,\pi_2\in \End_0(\cM)$, we denote the intertwiner space 
$(\pi_1,\pi_2):=\{T\in \cM\mid T\pi_1(x)=\pi_2(x)T \mbox{ for all } x\in \cM\}$.

Let $\rho\in \End_0(\cM)$, $\sigma\in \End_0(\cP)$ be endomorphisms 
such that $\cN=\rho(\cM)$, $\cQ=\sigma(\cP)$.
Fix standard isometries 
$R_\rho\in (\id, \bar{\rho}\rho )$,
$\bar{R}_\rho\in (\id, \rho\bar{\rho})$, 
$R_\sigma\in (\id, \bar{\sigma}\sigma )$ 
and $\bar{R}_\sigma\in (\id, \sigma\bar{\sigma})$ 
such that 
\[
 R^*_\rho\bar{\rho}(\bar{R}_\rho)=\bar{R}^*_\rho\rho(R_\rho)=\frac{1}{d(\rho)},\,\,
 R^*_\sigma\bar{\sigma}(\bar{R}_\sigma)=\bar{R}^*_\sigma{\sigma}(R_\sigma)=\frac{1}{d(\sigma)}.
\]

Set 
\[
\cM\supset \rho(\cM)\supset \rho \bar{\rho}(\cM)\supset 
\rho\bar{\rho}\rho(\cM)\supset 
\cdots =
\cM\supset \cN_1\supset \cN_2\supset \cN_3\cdots. 
     \]
\[
\cP\supset \sigma(\cP)\supset \sigma \bar{\sigma}(\cP)\supset 
\sigma\bar{\sigma}\sigma(\cP)\supset 
\cdots =
\cP\supset \cQ_1\supset \cQ_2\supset \cQ_3\cdots. 
     \]

Let $\gamma_\rho=\rho\bar{\rho}$, $\gamma_\sigma=\sigma \bar{\sigma}$ be
canonical endomorphisms.

Jones projections $e_n^\rho\in \cN_n$ ($n=0,1,2,\cdots$) are given as follows;
\[
 e_{2k}^\rho=\gamma_\rho^k(\bar{R}_\rho \bar{R}_\rho^*),\,\,\, 
 e_{2k+1}^\rho=\gamma_\rho^k\rho({R}_\rho {R}_\rho^*).
\]

Assume there exists an
isomorphism $\alpha\colon\rho(\cM)\subset\cM \rightarrow \sigma(\cP)\subset \cP$. 
In the definition of the Loi invariant, there is freedom of choice of a
unitary $u$. In what follows, we will fix the choice of a unitary by using standard
isometries. 

Since $\alpha\rho(\cM)=\sigma(\cP)$, 
there exists an isomorphism $\beta\colon\cM \rightarrow \cP$ such that
$\alpha\rho=\sigma\beta$. Thus $[\bar{\sigma}\alpha]=[\beta\bar{\rho}]$
holds as sectors,
and hence there exists $u\in U(\cP) $ such that $\Ad u\circ  \bar{\sigma}\alpha=
\beta\bar{\rho}$. Such unitary  $u$ is not uniquely
determined. However 
we can easily see 
$\sigma(u)^*\alpha(\bar{R}_\rho)\in (\id,\sigma\bar{\sigma})$,
$u^*\beta({R}_\rho)\in (\id,\bar{\sigma}\sigma)$, and these
isometries are also  standard solutions for $\sigma$,
$\bar{\sigma}$. So we can choose $u\in \cP$ in such a way
$\bar{R}_\sigma=\sigma(u)^*\alpha(\bar{R}_\rho)\in (\id,\sigma\bar{\sigma})$,
$R_{\sigma}=u^*\beta(R_\rho)\in (\id,\bar{\sigma}\sigma)$
\cite[Proposition 2.2.15]{Nes-Tu}.
In what follows, we fix these standard isometries and a unitary $u$ as above.

Define $v^{(0)}:=1$ and 
$v^{(k+1)}:=v^{(k)}(\sigma \bar{\sigma})^k\sigma(u)
(=v^{(k)}\gamma^k_\sigma(v^{(1)}))$.
Then $v^{(k-1)*}v^{(k)}\in \cQ_{2k-1}$, and 
$\alpha\circ (\rho\bar{\rho})^k=\Ad v^{(k)}\circ (\sigma
\bar{\sigma})^k\circ \alpha$ hold. 

Let $\alpha^{(k)}:=\Ad v^{(k)*}\circ \alpha$. Then 
$\alpha^{(k)}(e_l^\rho)=e_l^\rho$, $0\leq l\leq 2k-1$,
and  
$\alpha^{(k)}$sends $\cM\supset \cN_1\supset \cdots \supset \cN_{2k}\supset
\cN_{2k+1}$
to 
$\cP\supset \cQ_1\supset \cdots \supset \cQ_{2k}\supset
\cQ_{2k+1}$.
Thus $\alpha^{(k)}(\cN_l'\cap \cM)=\cQ_l'\cap \cP$ for all $0\leq l\leq 2k+1$.
Since $v^{(k-1)*}v^{(k)}\in \cQ_{2k-1}$, we have
$\alpha^{(k)}|_{\cN_{2k-1}'\cap \cM}=
\alpha^{(k-1)}|_{\cN_{2k-1}'\cap \cM}$. 
Thus $\{\alpha^{(k)}|_{\cN_{2k+1}'\cap \cM}\}_{k}$ gives a Loi invariant
for $\alpha$.

The Loi invariant $\Phi(\alpha)$ gives 
an isomorphism of standard invariants 
of $\rho(\cM)\subset \cM$ and $\sigma(\cP)\subset \cP$. 
Here we clarify the meaning of ``an isomorphism of standard invariants''
\cite[pp. 285]{Loi-auto}.

Let  $A_{n}^\rho=\cN_n'\cap \cM$, $B_{n}^\rho=\cN_{n}'\cap \cN$.
These relative commutants are described in terms of intertwiner spaces as follows;
\[
A_{2k}^\rho= (\gamma_\rho^{k},\gamma_\rho^{k}), \,\,\,
A_{2k+1}^\rho=
(\gamma_\rho^{k}{\rho},\gamma_\rho^{k}{\rho}),
\]
\[
B_{2k}^\rho=\rho\Bigl(
((\bar{\rho}\rho)^{k-1}\bar{\rho},(\bar{\rho}\rho)^{k-1}\bar{\rho})\Bigr), \,\,\,
B_{2k+1}^\rho=
\rho\Bigl(((\bar{\rho}\rho)^k,(\bar{\rho}\rho)^k)\Bigr).\]

We have a canonical inclusion 
\[
 \begin{array}{ccc}
  B_n^\rho&\subset  &B_{n+1}^\rho \\
  \cap& &\cap \\
  A_n^\rho &\subset  &A_{n+1}^\rho 
 \end{array}.
\]

Let $E_n^\rho$ be a conditional expectation 
$E_n^\rho\colon A_n^\rho\rightarrow A_{n-1}^\rho$  defined as follows.
\[
 E_{2k}(T)=\gamma_\rho^{k-1}\rho(R_\rho^*)
T\gamma_\rho^{k-1}\rho(R_\rho), \,\,
 E_{2k+1}(T)=\gamma_\rho^{k}(\bar{R}_\rho^*)
T\gamma_\rho^{k}(\bar{R}_\rho)
\]
It is easy to see $E_n^\rho \mid_{B_n^\rho}$ is a conditional expectation
$B_n^\rho\rightarrow B_{n-1}^\rho$.

Let 
$F^\rho_n(T):=\rho\bigl({R}_{\rho}^*\bar{\rho}(T){R}_{\rho}\bigr)$, $T\in
A_n^\rho$.  
Then $F_n^\rho$ is a conditional expectation from $A_n^\rho$ onto $B_n^\rho$.

Via these expectations, 
\[
 \begin{array}{rcl}
B_n^\rho&\overset{E_{n+1}^\rho}{\subset}  &B_{n+1}^\rho \\
\mbox{\raise3pt\hbox{${}_{F_{n}^\rho}$}}\cap\,\,& &
\,\cap\mbox{\raise3pt\hbox{${}_{F_{n+1}^\rho}$}} \\
  A_n^\rho &\underset{E_{n+1}^\rho}{\subset}  &A_{n+1}^\rho 
 \end{array}
\]
forms a commuting square.

The standard invariant of $\rho(\cM)\subset \cM$ is given by the
following nest of finite dimensional algebras
\[
\begin{array}{ccccccccc}
 B_1^\rho&\subset  &B_2^\rho
  &\subset  & B_3^\rho
  &\subset  & B_4^\rho&\cdots \\
 \cap & &\cap & &\cap & &\cap & \\
  A_1^\rho&\subset  &A_2^\rho
  &\subset  & A_3^\rho
  &\subset  & A_4^\rho &\cdots 
\end{array} 
\]
together with conditional expectations and Jones projections.

An isomorphism of  standard invariants from $\rho(\cM)\subset \cM$ to
 $\sigma(\cP)\subset \cP$
 is a family of maps
$\{\alpha_k\}_{k=1}^\infty$, such that 
\begin{enumerate}
\setlength{\itemsep}{-5pt}
 \item $\alpha_k$ is an  isomorphism from $A_k^\rho\supset B_k^\rho$ onto
       $A_k^\sigma\supset B_k^\sigma$ such that $F_k^\sigma\alpha_k=\alpha_kF_k^\rho$.  
 \item $\alpha_{k}\mid_{A_{k-1}^\rho}=\alpha_{k-1} $, 
$E_k^\sigma\alpha_k=\alpha_{k-1}E_k^\rho$. 
 \item $\alpha_{k+2}(e_k^\rho)=e_k^\sigma$. (Note $e_k^\rho\in A_{k+2}$.)
\end{enumerate}
Note that 
$\alpha_k$ is an isomorphism of  commuting squares 
\[
 \begin{array}{ccc}
  B_{k-1}^\rho&\subset  &B_k^\rho \\
	 \cap& &\cap \\
 A_{k-1}^\rho&\subset  &A_k^\rho
	\end{array}
\longrightarrow 
 \begin{array}{ccc}
  B_{k-1}^\sigma&\subset  &B_k^\sigma \\
		 \cap& &\cap \\
A_{k-1}^\sigma&\subset  &A_k^\sigma
\end{array}.
\]

Here we assume $\cN\subset \cM$ 
and $\cQ\subset \cP$ are strongly
amenable subfactors of type III$_1$ in the sense of Popa 
\cite{Po-amen},
\cite{Po-tani}, \cite{Po-CBMS},  with identical  type II principal graph
and type III graph. Then these subfactors are classified by their standard
invariants \cite{Po-CBMS} (also see \cite{M-III1}.) We state this
classification more precisely.
There exists a subfactor of type II$_1$ $\cB\subset \cA$, 
  such that $\cB\subset \cN$, $\cA\subset \cM$, and a 
tunnel 
\[
\cA\supset \cB\supset \cB_1\supset \cdots,\,\,
\]
such that $\cB\subset \cA=\cB^{\mathrm{st}}\subset \cA^{\mathrm{st}}$, 
and 
$\cN\subset \cM=(\cB^{\mathrm{st}}\subset
\cA^{\mathrm{st}})\otimes \cR_\infty$.
Here
$\cA^{\mathrm{st}}=\bigvee_{k}(\cB_k'\cap \cA)$, 
$\cB^{\mathrm{st}}=\bigvee_{k}(\cB_k'\cap \cB)$, and $\cR_\infty$ is the
unique injective factor of type III$_1$ \cite{Co-III1}, \cite{Ha-III1}.

Let $\cT\subset \cS$ be a subfactor of type II$_1$ which has a same
property for $\cQ\subset \cP$ as above.
If we have an isomorphism $\theta_0$ on the standard invariants of $\cN\subset
\cM$ and $\cQ\subset \cP$, then it extends to an isomorphism from
$\cB^{\mathrm{st}}\subset \cA^{\mathrm{st}}$ onto $\cT^{\mathrm{st}}\subset \cS^{\mathrm{st}}
$, and hence to that of 
$\cN\subset \cM$ onto $\cQ\subset \cP$. Moreover, the Loi invariant of
this isomorphism coincides with $\theta_0$.
We summarize this as follows.
\begin{thm}\label{thm:isomorphism}
 Let $\theta_0$ be an isomorphism of standard invariants of $\cN\subset
 \cM$ and 
$\cQ\subset \cP$. Then there exists an isomorphism $\theta\colon\cN\subset
 \cM\rightarrow \cQ\subset \cP$ whose Loi
 invariant is $\theta_0$. 
\end{thm}

\noindent
\textbf{Remark.} The above $\theta$ does not necessary preserve given
tunnel for  subfactors. 
We can perturb $\theta$ so that it preserves any
finite length of a tunnel.

\section{Equivalence and cocycle conjugacy}\label{sec:equi}

In \cite{Rob-act}, Roberts defined an action of a group dual on a von
Neumann algebra as a monoidal functor from a representation category of
a group to a endomorphism category of a von Neumann algebra. By
generalizing his definition, we introduce the notion of a Roberts action
of a C$^*$-tensor category. 

\begin{df}\label{df:Roberts action}
\upshape
Let $\mathcal{C}$ be a C$^*$-tensor category, and $\cM$ a factor. \\
{\upshape (1)} A Roberts action $\rho$ of $\mathcal{C}$ on $\cM$ is a
 C$^*$-tensor functor $\rho\colon\xi\in \mathcal{C}\mapsto \rho_\xi\in
 \End_0(\cM)$ such  that $\rho_{\xi\otimes \eta}=\rho_\xi\rho_\eta$. \\
{\upshape (2)} A Roberts action $\rho$ of $\mathcal{C}$ is said to be free if $\rho$ is a
 fully faithful 
 functor, i.e., $\rho\colon T\in\Mor(\xi,\eta)\mapsto \rho_T\in (\rho_\xi,\rho_\eta)$
 is bijective for any $\xi,\eta\in \mathcal{C}$. \\ 
{\upshape (3)} A Roberts action $\rho$ of $\mathcal{C}$ is said to be modularly
 free if it is free and 
$\tilde{\rho}_\xi\in \End_0(\tM)$ is not a modular endomorphism
in the sense of {\upshape \cite[Definition 3.1]{Iz-can2}} for
 $\xi\not\cong 1_\mathcal{C}$. 
Here $\tM$ is the core for $\cM$, and 
$\tilde{\rho}_\xi$ is the canonical extension of $\rho_\xi$ in the
 sense of 
{\upshape \cite[Theorem 2.4]{Iz-can2}}.
\end{df}

\noindent
\textbf{Remark.}
(1) In the rest of this paper,
we use a same letter
$T\in (\rho_\xi,\rho_\eta)$  for an intertwiner $T\in \Mor(\xi,\eta)$
to simplify notation. \\
(2) When $\cM$ is a factor of type III$_1$, 
the modular freeness of $\rho$ is equivalent to the 
full faithfulness of a functor $\xi \in \mathcal{C}\mapsto
\tilde{\rho}_\xi\in \End_0(\tM)$. 
In general, 
the modular freeness of $\rho$ is equivalent to  
$(\tilde{\rho}_\xi,\tilde{\rho}_\eta)=Z(\tM)\otimes
(\rho_\xi,\rho_\eta)$ for $\xi,\eta\in \mathcal{C}$. 
(Here we have
$Z(\tM), (\rho_\xi,\rho_\eta)\subset
(\tilde{\rho}_\xi,\tilde{\rho}_\eta)$ 
and $Z(\tM)\cap (\rho_\xi,\rho_\eta)=\mathbb{C}$.
Hence we can regard $Z(\tM)\otimes (\rho_\xi,\rho_\eta)\subset
(\tilde{\rho}_\xi,\tilde{\rho}_\eta)$.) \\
(3) When $\cM$ is injective, the modular freeness is equivalent to
the central freeness \cite[Theorem 1]{KwST}, \cite[Theorem 4.12]{MaTo-endo}.

\medskip

Let $\mathcal{C}$, $\mathcal{D}$ be  unitarily equivalent 
C$^*$-tensor categories, and 
$(F,L)$ be a tensor equivalence functor of $\mathcal{C}$ and
$\mathcal{D}$. Namely, $F$ is a functor $F\colon\mathcal{C}\rightarrow
\mathcal{D}$, which gives surjective isomorphisms on intertwiner spaces, and $L$ is
a natural unitary equivalence of tensor product
$L(\xi,\eta)\in (F(\xi)\otimes F(\eta), F(\xi\otimes \eta))$.
We assume $F(1_{\mathcal{C}})=1_{\mathcal{D}}$, and $L(1_{\mathcal{C}},\xi)=L(\xi,1_{\mathcal{C}})=1$.
In the following, we omit $\otimes$, and 
denote $\xi\otimes \eta$ by $\xi\eta$ for
simplicity.

Let $\rho$ (resp. $\sigma$) be 
a modularly free Roberts action of $\mathcal{C}$ on $\cM$ (resp. $\mathcal{D}$ on $\cP$),
where $\cM$ and $\cP$ are factors of type III$_1$. 
Such actions exist by the result of Hayashi-Yamagami \cite[Theorem 7.6]{HaYa-amenable} together with a simple
trick by making tensor product with an injective factor of type III$_1$,
if $\mathcal{C}$ and $\mathcal{D}$ are amenable and factors are injective.

By the remark (1) after Definition \ref{df:Roberts action},  an element of $(\sigma_{F(\xi)},\sigma_{F(\eta)})$
can be  expressed by $F(T)$, $T\in(\rho_\xi,\rho_\eta)$, and  $L(\xi,\eta)$
can be regarded as  a unitary in
$(\sigma_{F(\xi)}\sigma_{F(\eta)},\sigma_{F(\xi\eta)})$.

We extend $L(\xi,\eta)$ to $L(\xi_1,\cdots ,\xi_n)\in 
(\sigma_{F(\xi_1)}\cdots\sigma_{F(\xi_n)},\sigma_{F(\xi_1\cdots
\xi_n)})$ in the canonical way. We
use the notation $L(\xi)=1$ for  $\xi\in \mathcal{C}$. For $T\in
(\rho_{\xi_1}\cdots\rho_{\xi_n},
\rho_{\eta_1}\cdots\rho_{\eta_m})$,
let 
\[
 \theta(T):=L(\eta_1,\cdots ,\eta_m)^*F(T)L(\xi_1,\cdots ,\xi_n)
\in 
(\sigma_{F(\xi_1)}\cdots\sigma_{F(\xi_n)},
\sigma_{F(\eta_1)}\cdots\sigma_{F(\eta_m)}).\]
Although we should use the notation like
$\theta_{\xi_1,\cdots,\xi_n}^{\eta_1,\cdots ,\eta_m}(T)$ for completeness,
we simply denote by $\theta(T)$.

\begin{lem}\label{lem:Fiso}
We have the followings. \\
$(1)$ We have $\theta(ST)=\theta(S)\theta(T)$ and
 $\theta(T)^*=\theta(T^*)$ for $T\in (\rho_{\xi_1}\cdots\rho_{\xi_n},
 \rho_{\eta_1}\cdots\rho_{\eta_m})$,
 $S\in (\rho_{\eta_1}\cdots\rho_{\eta_m},\rho_{\zeta_1}\cdots
 \rho_{\zeta_l})$. \\
 $(2)$ $\theta(T)$ is well-defined in the following sense.
If we regard $T\in (\rho_{\xi_1}\cdots\rho_{\xi_n},
\rho_{\eta_1}\cdots\rho_{\eta_m})$ as an element  in $
(\rho_{\xi_1}\cdots\rho_{\xi_n}\rho_\zeta,
\rho_{\eta_1}\cdots\rho_{\eta_m}\rho_\zeta)$, then we have
\[
 L(\eta_1,\cdots ,\eta_m)^*F(T)L(\xi_1,\cdots ,\xi_n)=
L(\eta_1,\cdots ,\eta_m,\zeta)^*F(T)L(\xi_1,\cdots ,\xi_n,\zeta).
\]
$(3)$ $\sigma_{F(\zeta)}\theta(T)=\theta(\rho_{\zeta}(T))$ 
for $T\in (\rho_{\xi_1}\cdots\rho_{\xi_n},
 \rho_{\eta_1}\cdots\rho_{\eta_m})$. \\
\end{lem}

\noindent
\textbf{Proof.}
(1) It is easy to see  $\theta(T)^*=\theta(T^*)$. 
We can verify   $\theta(TS)=\theta(T)\theta(S)$  as follows.
\begin{align*}
\theta(S)\theta(T)
 &=L(\zeta_1,\cdots \zeta_l)^*F(S)L(\eta_1,\cdots, \eta_m)
L(\eta_1,\cdots, \eta_m)^*F(T)L(\xi_1,\cdots, \xi_n) \\
&=L(\zeta_1,\cdots \zeta_l)^*F(S)F(T)L(\xi_1,\cdots, \xi_n) \\
&=\theta(ST). 
\end{align*}
(2) Note $L(\xi_1,\cdots, \xi_n,\zeta)=L((\xi_1\cdots\xi_n),\zeta)L(\xi_1,\cdots,\xi_n)$.
Then 
\begin{align*}
\lefteqn{L(\eta_1,\cdots ,\eta_m,\zeta)^*F(T)L(\xi_1,\cdots
 ,\xi_n,\zeta)} \\ &=
L(\eta_1,\cdots ,\eta_m)^*L((\eta_1\cdots\eta_m),\zeta)^*
F(T)L((\xi_1\cdots\xi_n),\zeta)L(\xi_1,\cdots ,\xi_n) \\ 
&=
L(\eta_1,\cdots ,\eta_m)^*L((\eta_1\cdots\eta_m),\zeta)^*L((\eta_1\cdots\eta_m),\zeta)
F(T)L(\xi_1,\cdots ,\xi_n)\,\, \mbox{(by naturality of $L$)}
 \\ 
&=
L(\eta_1,\cdots ,\eta_m)^*F(T)L(\xi_1,\cdots ,\xi_n).
\end{align*}
(3)
Note $L(\zeta,\xi_1,\cdots, \xi_n)=L(\zeta,(\xi_1\cdots\xi_n))\sigma_{F(\zeta)}(L(\xi_1,\cdots,\xi_n))$.
Then 
\begin{align*}
\lefteqn{\theta(\rho_\zeta(T))} \\
&=
L(\zeta,\eta_1,\cdots ,\eta_m)^*F(\rho_\zeta(T))L(\zeta,\xi_1,\cdots
 ,\xi_n) \\ 
&=
 \sigma_{F(\zeta)}(L(\eta_1,\cdots,\eta_m))^*
L(\zeta ,(\eta_1\cdots\eta_m))^*
F(\rho_\zeta(T))
L(\zeta, (\xi_1\cdots\xi_n)) \sigma_{F(\zeta)}\left(L(\xi_1,\cdots,\xi_n)\right)\\ 
&=
 \sigma_{F(\zeta)}(L(\eta_1,\cdots,\eta_m))^*
L(\zeta ,(\eta_1\cdots\eta_m))^* L(\zeta, (\eta_1\cdots\eta_m)) 
\sigma_{F(\zeta)}(F(T))
\sigma_{F(\zeta)}\left(L(\xi_1,\cdots,\xi_n)\right) \\ 
&=
 \sigma_{F(\zeta)}(L(\eta_1,\cdots,\eta_m))^*
\sigma_{F(\zeta)}(F(T))
\sigma_{F(\zeta)}L(\xi_1,\cdots,\xi_n) \\
&=\sigma_{F(\zeta)}(\theta(T)).
\end{align*}
Here the third equality is due to the naturality of $L$. \hfill$\Box$

\medskip

We present the definition of cocycle conjugacy of two Roberts actions of
C$^*$-tensor category.
\begin{df}\label{df:equivalence}
\upshape
 Let $\mathcal{C}$ and  $\mathcal{D}$ be unitarily equivalent C$^*$-tensor categories with
 tensor equivalence $(F,L)$. Let $\rho$ (resp. $\sigma$) be a Roberts action
 of $\mathcal{C}$ (resp. $\mathcal{D}$) on a factor $\cM$ (resp. $\cP$).
We say $\rho$ and $\sigma$ are cocycle conjugate if 
there exist an isomorphism $\alpha\colon\cM\rightarrow\cP$ and 
a unitary $\mathcal{E}(\xi)\in \cP$ for every $\xi \in
 \mathcal{C}$
 such that \\
$\mathrm{(i)}$ $\mathcal{E}(1_{\mathcal{C}})=1$, \\
$\mathrm{(ii)}$ $\mathcal{E}(\xi)\in
 (\sigma_{F(\xi)}\alpha,\alpha\rho_\xi)$, \\
$\mathrm{(iii)}$
 $\mathcal{E}(\xi)\sigma_{F(\xi)}(\mathcal{E}(\eta))=\mathcal{E}(\xi\eta)L(\xi,\eta)$, \\
$\mathrm{(iv)}$ $\alpha(T)\mathcal{E}(\xi)=\mathcal{E}(\eta)F(T)$ for $T\in
 \Mor(\xi,\eta)$. In particular,
 $\mathcal{E}(\xi)=\sum \alpha(T_i)\mathcal{E}(\xi_i)F(T_i)^*$ holds for
 a (not necessary irreducible)
 decomposition $\rho_\xi(x)=\sum_i T_i\rho_{\xi_i}(x)T_i^*$. \\
\end{df}
\textbf{Remark.} In Definition \ref{df:equivalence}, we do not assume
the simplicity of objects. \\

Assume $\mathcal{C}$ and $\mathcal{D}$ are finitely generated.
Let $\mathcal{C}_0\subset \mathrm{Irr}(\mathcal{C})$ be a (finite)
generator of $\mathcal{C}$. 
%and $\mathcal{D}_0:=F(\mathcal{C}_0)$. 
We assume that $1_{\mathcal{C}}\in \mathcal{C}_0$ and $\mathcal{C}_0$ is closed under conjugation.

Let $\pi:=\bigoplus_{\xi\in \mathcal{C}_0}\xi$, and 
consider subfactors $\rho_{\pi}(\cM)\subset \cM$ and $\sigma_{F(\pi)}(\cP)\subset\cP$.
(Note $\rho_\pi$ is self-conjugate.)
Set $\gamma_\rho:=\rho_\pi\rho_\pi$, $\gamma_{\sigma}:=\sigma_{F(\pi)}\sigma_{F(\pi)}$.

We fix standard
isometries $R_{\rho}\in(\id, \bar{\rho}_\pi\rho_\pi)$,
$\bar{R}_{\rho}\in(\id, {\rho}_\pi\bar{\rho}_\pi)$,
${R}_{\sigma}\in(\id, \bar{\sigma}_{F(\pi)}{\sigma}_{F\pi})$,
$\bar{R}_{\sigma}\in(\id, {\sigma}_{F(\pi)}\bar{\sigma}_{F(\pi)})$
so that $\theta(R_{\rho})=R_\sigma$,
$\theta(\bar{R}_\rho)=\bar{R}_\sigma$.

Standard invariants of these subfactors depends only on 
$\mathcal{C}$, $\mathcal{D}$ and not on actions $\rho,\sigma$.
In fact, by Lemma \ref{lem:Fiso}
we can see that  
an isomorphism between standard invariants is given by 
\[
\theta\colon(\rho_\pi^n,\rho_\pi^n)\rightarrow (\sigma_{F(\pi)}^n,\sigma_{F(\pi)}^n), 
\]
and hence by an equivalence
$(F,L)$. Since $\rho$ and $\sigma$ are modularly free, 
these subfactors have same type II principal graphs
and type III graphs by \cite[Theorem 3.5]{Iz-graph}.

Let us assume $\mathcal{C}$ and $\mathcal{D}$ are strongly amenable,
i.e.,  fusion algebras $\mathbb{C}[\mathrm{Irr}(\mathcal{C})]$, 
$\mathbb{C}[\mathrm{Irr}(\mathcal{D})]$ are strongly amenable in the
sense of \cite[Definition 6.4]{Hi-Iz}. 
By \cite[Theorem 4.8]{Hi-Iz}, 
this is equivalent to the strong amenability of 
standard invariants of $\rho_\pi(\cM)\subset \cM$ and
$\sigma_{F(\pi)}(\cP)\subset \cP$ 
in the sense of Popa \cite{Po-amen}, \cite{Po-tani}, \cite{Po-CBMS}.

We state the main theorem of this paper.

\begin{thm}\label{thm:equivalence}
 Let $\mathcal{C}$ and  $\mathcal{D}$ be unitarily equivalent finitely
 generated strongly amenable C$^*$-tensor categories with
 tensor equivalence $(F,L)$. Let $\rho$ (resp. $\sigma$) be a  modularly free Roberts action
 of $\mathcal{C}$ (resp. $\mathcal{D}$) on an injective factor of type III$_1$
 $\cM$ (resp. $\cP$).  Then $\rho$ and $\sigma$ are cocycle conjugate.
\end{thm}

The following corollary immediately follows if we put
$\mathcal{C}=\mathcal{D}$, $\cM=\cP$, and $(F,L)$ as an identity functor.

\begin{cor}\label{cor:equivalence}
 Let $\mathcal{C}$
  be a finitely  generated strongly amenable C$^*$-tensor category. 
 Let $\rho$ and $\sigma$ be  modularly free Roberts actions 
of $\mathcal{C}$ on the injective factor $\cM$ of type III$_1$. 
Then there exist $\alpha \in \Aut(\cM)$ and a unitary
 $\mathcal{E}(\xi)$ for every $\xi\in \cM$ such that \\
$\mathrm{(i)}$ $\Ad \mathcal{E}(\xi)\circ \sigma_\xi
 =\alpha\circ\rho_\xi\circ\alpha^{-1}$ for every $\xi\in \mathcal{C}$. \\
$\mathrm{(ii)}$
 $\mathcal{E}(\xi)\sigma_{F(\xi)}(\mathcal{E}(\eta))=\mathcal{E}(\xi\eta)$. \\
$\mathrm{(iii)}$ $\alpha(T)\mathcal{E}(\xi)=\mathcal{E}(\eta)F(T)$ for $T\in
 \Mor(\xi,\eta)$. 
\end{cor}

\noindent
\textbf{Proof of Theorem \ref{thm:equivalence}.} 
The proof is a simple modification of proof of \cite[Theorem 2.1]{M-ExtIJM}.

By Theorem \ref{thm:isomorphism}, there exists an isomorphism
$\alpha\colon\rho_\pi(\cM)\subset \cM\rightarrow \sigma_{F(\pi)}(\cP)\subset
\cP$ whose Loi invariant  is $\theta$. Then we will choose $u\in \cP$
so that $\bar{R}_\sigma=\sigma(u^*)\alpha(\bar{R}_\rho)$ as explained in
\S\ref{sec:half}.

Let us  summarize  facts  which are used in the
proof.
\[
v^{(1)}=\sigma(u),\,\,
v^{(n+1)}=v^{(n)}\gamma_\sigma^n(v^{(1)})=v^{(1)}\gamma_\sigma(v^{(n)}), \,\,
\]
\[
\alpha\circ \gamma_\rho^{n}=\Ad v^{(n)}\circ\gamma_\sigma^n\circ\alpha, \,\,
 \alpha^{(n)}=\Ad v^{(n)*}\circ \alpha, 
\]
\[
 \alpha^{(n)}=\theta \mbox{ on } (\gamma_\rho^n,\gamma_\rho^n),
\]
\[
 \theta(\bar{R}_\rho)=\bar{R}_\sigma=\sigma(u^*)\alpha(\bar{R}_\rho)=
v^{(1)*}\alpha(\bar{R}_\rho).
\]

Let $\xi \in \mathcal{C}$. Since $\mathcal{C}_0$ generates
$\mathcal{C}$, there exist $n\in \mathbb{N}$ and an isometry
$T\in (\rho_\xi, \gamma_\rho^{n})$. 

Define 
$W_T=\alpha(T^*)v^{(n)}\theta(T)$. We will verify that $W_T$ is a
desired unitary $\mathcal{E}(\xi)$.

Our first task is to 
show that $W_T$ is a unitary and independent from the choice of $T$
and $n$. 

First, we show $W_T=W_S$ for 
$T, S\in (\rho_\xi, \gamma^n_\rho)$. We compute $W_T^*W_S=1$. We have
\[
 W_T^*W_S=
\theta(T^*) v^{(n)*} \alpha(T)
\alpha(S^*)v^{(n)}\theta(S) =
\theta(T^*)\alpha^{(n)}(TS^*)\theta(S).
\]
Here  $TS^*\in (\gamma_\rho^n,\gamma_\rho^n)$, and 
$\alpha^{(n)}$ on $(\gamma_\rho^n, \gamma_\rho^n)$ is
given by $\theta$. Hence
\[
 \theta(T^*)\alpha^{(n)}(TS^*)\theta(S)=\theta(T^*)\theta(TS^*)\theta(S)=1
\]
holds. Similarly we get $W_SW_T^*=1$. 
This shows that $W_T$ is a unitary (put $T=S$), and $W_T=W_S$.

Next we show that $W_T$ does not depend on $n$. The following proof is
suggested by R. Tomatsu, which is much simpler than author's original
proof. 
Let $T\in (\rho_\xi,\gamma_\rho^n)$ be an isometry. Then
$\bar{R}_\rho T\in (\rho_\xi, \gamma_\rho^{n+1})$. 
\begin{align*}
 W_{\bar{R}_\rho T}&=\alpha(T^*\bar{R}_\rho^*)v^{(n+1)}
\theta(\bar{R}_\rho T) \\ 
&=
\alpha(T^* \bar{R}_\rho^*)v^{(1)}\gamma_{\sigma}(v^{(n)})\theta(\bar{R}_\rho T) \\
&=
\alpha(T^*)\bar{R}_\sigma^*
\gamma_{\sigma}(v^{(n)})\theta(\bar{R}_\rho)\theta( T) \\
&=
\alpha(T^*)
v^{(n)}\bar{R}_\sigma^*\bar{R}_{\sigma}\theta(T)   \\ 
&=
\alpha(T^*)v^{(n)} \theta(T) \\
&=W_T.
\end{align*}

Thus $\mathcal{E}(\xi):=W_T$ is well-defined.

We will show that $\mathcal{E}(\xi)$ satisfies the condition (i)--(iv) in
Definition \ref{df:equivalence}. 

(i)
Since $\bar{R}_\rho\in (\id,\gamma_\rho)$, 
\[
 \mathcal{E}(1_{\mathcal{C}})=W_{\bar{R}_\rho}=\alpha(\bar{R}_\rho^*)v^{(1)}\theta(\bar{R}_\rho)=
\bar{R}_\sigma^*\bar{R}_\sigma=1.
\]

(ii) Choose $n\in \mathbb{N}$ and an isometry $T\in
(\rho_\xi,\gamma_\rho^n)$. 
Note $\theta(T)\in (\sigma_{F(\xi)},\gamma_\sigma^n)$. 
Then
$\mathcal{E}(\xi)=W_T\in (\sigma_{F(\xi)}\alpha,\alpha\rho_\xi)$ is verified as follows;
\begin{align*}
 \alpha\rho_\xi(x)W_T&=
 \alpha\rho_\xi(x)\alpha(T^*)v^{(n)}\theta(T) \\
&=\alpha(T^*)
 \alpha\gamma_\rho^n(x)
v^{(n)}\theta(T) \\
&=\alpha(T^*)v^{(n)}
 \gamma_\sigma^n\alpha(x)\theta(T) \\
&=\alpha(T^*)v^{(n)}\theta(T)
 \sigma_{F(\xi)}\alpha(x) \\
&=W_T
 \sigma_{F(\xi)}\alpha(x).
\end{align*}

To prove (iii) and (iv), we claim the following.\\
\textbf{Claim.} We have
\[
\mathcal{E}(\xi)\sigma_{F(\xi)}\left(\mathcal{E}(\eta)\right)\theta(T)=
\alpha(T)\mathcal{E}(\zeta) 
\]
for $\xi,\eta,\zeta\in \mathcal{C}$ and 
an isometry $T\in (\rho_\zeta,\rho_\xi\rho_\eta)$.

Take $n,m\in \mathbb{N}$ and isometries $S_\xi\in (\rho_\xi,\gamma_\rho^n)$,
$S_\eta\in (\rho_\eta,\gamma^m_\rho)$. 
Then $S_\zeta:=S_\xi\rho_\xi(S_\eta)T\in (\rho_\zeta,\gamma^{n+m}_\rho)$
is an isometry, and $S_\zeta S_\zeta^*\in (\gamma_\rho^{n+m},\gamma_{\rho}^{n+m})$.
So we have $\alpha^{(n+m)}(S_\zeta S_\zeta^*)=\theta(S_\zeta S_\zeta^*)$.

We show $W_{S_\xi}\sigma_{F(\xi)}(W_{S_\eta})\theta(T)=\alpha(T)W_{S_\zeta}$;
\begin{align*}
\lefteqn{ W_{S_\xi}\sigma_{F(\xi)}(W_{S_\eta})\theta(T)} \\
&=W_{S_\xi}
\sigma_{F(\xi)}\alpha(S_\eta^*)\sigma_{F(\xi)}(v^{(m)}) \sigma_{F(\xi)}(\theta(S_\eta))\theta(T) \\
&=
\alpha\rho_\xi(S_\eta^*)W_{S_\xi}
\sigma_{F(\xi)}(v^{(m)}) \sigma_{F(\xi)}(\theta(S_\eta))\theta(T)  \\
&=
\alpha\rho_\xi(S_\eta^*)
\alpha(S_\xi^*)v^{(n)} \theta(S_\xi)
\sigma_{F(\xi)}(v^{(m)}) \sigma_{F(\xi)}(\theta(S_\eta))\theta(T) \\
&=
\alpha\rho_\xi(S_\eta^*)
\alpha(S_\xi^*)v^{(n)}
\gamma^n_\sigma(v^{(m)})
 \theta(S_\xi)\theta(\rho_\xi(S_\eta))\theta(T)\,\,\,\,\mbox{(by Lemma \ref{lem:Fiso}(2))}
 \\
&=
\alpha\left(\rho_\xi(S_\eta^*)S_\xi^*\right)v^{(n+m)}
\theta(S_\xi\rho_\xi(S_\eta)T)
  \\
&=v^{(n+m)} 
\alpha^{(n+m)}\left(\rho_\xi(S_\eta^*)S_\xi^*\right)
\theta(S_\zeta)  \\
&=v^{(n+m)}
\alpha^{(n+m)}\left(\rho_\xi(S_\eta^*)S_\xi^*\right)
\theta(S_\zeta S_\zeta^*S_\zeta)\,\,\,(\mbox{since }S_\zeta\mbox{ is an isometry})   \\
&=v^{(n+m)}
\alpha^{(n+m)}\left(\rho_\xi(S_\eta^*)S_\xi^*
S_\zeta S_\zeta^*
\right)
\theta(S_\zeta )
\\
&=v^{(n+m)}
\alpha^{(n+m)}\left(TS_\zeta^*\right)
\theta(S_\zeta )
\\
&=\alpha(T)
\alpha\left(S_\zeta^*\right)v^{(n+m)}
\theta(S_\zeta )
\\
&=\alpha(T)W_{S_\zeta}.
\end{align*}
Thus Claim has been verified. We continue the proof of Theorem \ref{thm:equivalence}.

(iii)   Put $\zeta=\xi\eta$, and $T=1\in
(\rho_{\xi\eta},\rho_\xi\rho_\eta)$ in Claim. 
Note $\theta(T)=L(\xi,\eta)^*F(1)L(\xi\eta)=
L(\xi,\eta)^*$ in this case. Then we immediately get (iii).

(iv) If $(\xi,\eta)\ne 0$, then there exists
$\zeta\in \mathcal{C}$ such that
$\zeta\prec\xi$, $\zeta\prec\eta$. Take isometries $T\in\Mor(\zeta,\xi)$ and 
$S\in\Mor(\zeta,\eta)$. Then $ST^*\in \Mor(\xi,\eta)$. We only have to
show (iv) for $ST^*$, since every element in $\Mor(\xi,\eta)$ is a  linear
span of intertwiners of such form.

By Claim, 
we have
\[
 \mathcal{E}(\xi)\theta(T)=\alpha(T)\mathcal{E}(\zeta), \,\,\,
 \mathcal{E}(\eta)\theta(S)=\alpha(S)\mathcal{E}(\zeta).
\]
(Put $\xi=1_\mathcal{C}$ or $\eta=1_\mathcal{C}$ in Claim.)
Then 
\begin{align*}
\alpha(ST^*)\mathcal{E}(\xi)=
\alpha(S) \mathcal{E}(\zeta)\theta(T)^*=
\mathcal{E}(\eta)\theta(ST^*)=
\mathcal{E}(\eta)F(ST^*).
\end{align*}
Note $\theta(ST^*)=L(\eta)^*F(ST^*)L(\xi)=F(ST^*)$.

\hfill $\Box$

\bigskip

\noindent
\textbf{Remark}. 
(1) Even when $\cM$ and $\cP$ are general AFD type III
factors, Theorem \ref{thm:equivalence} 
is true provided  
that $\rho_\xi$ and $\sigma_\xi$ have trivial Connes-Takesaki modules
in the sense of \cite{Iz-can2} thanks to the  classification results 
\cite[Theorem 6.1]{Loi-auto}, \cite[Theorem 4.2]{Win-30}
of subfactors of type III$_\lambda$, $\lambda\ne 1$. \\
(2) We can easily see that the statement in Claim  
is  equivalent to the conditions (iii), (iv) in Definition
\ref{df:equivalence}
under the condition $\mathcal{E}(1_\mathcal{C})=1$.

\bigskip

We present typical applications of our main theorem.

\begin{exam}
\upshape
When $\mathcal{C}$ is a finite C$^*$-tensor
category, a free Roberts action of $\mathcal{C}$ is 
 automatically  modularly free. Thus Theorem \ref{thm:equivalence} is
a generalization of \cite[Theorem 2.2]{Izumi-neargrp}.
\end{exam}

Next application is  
classification of group actions.

\begin{exam}
\upshape
 Let $\Gamma$ be a discrete  group. Then a Roberts action of
 $\mathrm{Vec}_{\Gamma}$ is nothing but a usual action of $\Gamma$.
If $\Gamma$ is a finitely generated strongly amenable group, then the main
 theorem implies that two centrally free actions of $\Gamma$ on the
 injective factor of type III$_1$ are cocycle conjugate. 
(Note that complete classification has been already obtained in \cite{KwST},
 \cite{KtST}, \cite{M-unif-Crelle} by different methods.)

It is worth mentioning the following remark.
If we consider only an isomorphism between locally trivial subfactors, 
then we  get only outer conjugacy of actions.
To obtain cocycle conjugacy, one must examine the Loi invariant carefully. 
\end{exam}

Other interesting application of the main theorem is Popa-Wasserman's
theory on classification of compact Lie group actions.
\begin{exam}[\cite{PoWa}]
\upshape
Let $\mathbb{G}$ be a compact Kac algebra such that $\mathrm{Rep}(\mathbb{G})$
 is strongly amenable, and finitely generated. Typical examples of such
 compact Kac algebras come from compact Lie groups \cite{Po-amen}, \cite{Hi-Iz}.
Then the main theorem implies the uniqueness
 of centrally free Roberts actions of $\mathrm{Rep}(G)$ on the injective factor of
 type III$_1$. Let $\pi$ be a generator of $\mathrm{Rep}(\mathbb{G})$,
 and set $\tilde{\pi}:=\id \oplus \pi \oplus \bar{\pi}$.
 Let $\alpha$ be a modularly free action of
 $\mathrm{Rep}(\mathbb{G})$ on an injective factor $\cM$ of type III$_1$ and  
$\hat{\alpha}$ a dual action of $\mathbb{G}$ on $\cN=\cM\rtimes_\alpha \hat{\mathbb{G}}$.
Consider a Wasserman subfactor $\cN^{\hat{\alpha}}\subset
(\cN\otimes M_{d\tilde{\pi}}(\mathbb{C}))^{\hat{\alpha}\otimes \Ad \tilde{\pi}}$. This
 is isomorphic to the dual inclusion of $\alpha_{\tilde{\pi}}(\cM)\subset
 \cM$. 
Thus main theorem
 implies $\alpha$ is unique up to cocycle conjugacy, and hence
 $\hat{\alpha}$ is unique up to conjugacy.

Note that centrally free actions of $\mathrm{Rep}(\mathbb{G})$ on the
 injective factor of type III$_1$ is classified completely in 
\cite[Theorem 7.16]{Mato-disKac-Memo} for a general coamenable Kac algebra
 $\mathbb{G}$ by a different method.
\end{exam}

We close this paper with the following comments.

As long as we apply classification results of subfactors, it seems to be
impossible to remove the assumption that tensor categories are finitely
generated. Thus  it is desirable to develop methods used in
\cite{MaTo-JFA}, \cite{Mato-disKac-Memo} to generalize Theorem
\ref{thm:equivalence} for arbitrary amenable C$^*$-tensor categories. \\

\noindent
\textbf{Problem.} For amenable C$^*$-tensor categories, show the
uniqueness of modularly free  Roberts actions on the injective factor of type
III$_1$.

\ifx\undefined\bysame
\newcommand{\bysame}{\leavevmode\hbox to3em{\hrulefill}\,}
\fi

\end{document}